\documentclass[conference]{IEEEtran}

\usepackage{graphicx}
\usepackage{amsmath, amsthm, amssymb, amsfonts}
\usepackage{algorithm, algorithmic, ifsym, subfigure, bm, empheq}

\usepackage{mathtools}

\newcommand{\defequal}{ \stackrel{\rm def}{=}  }

\newcommand{\figdep}[1]{}

\newcommand{\doFigure}[4]{
\ifnum\macInUse=1
  \begin{figure}[htbp]
    \begin{center}
      \vspace{.2in}
      \centerline {
        \epsfig{figure=#1,width=8cm}
      }
      \vspace{.2in}
      \caption{#2}
      \label{fig:#3}
    \end{center}
  \end{figure}
\else
\begin{figure}[htb]
  \centerline{\psfig{figure=#1.eps,width=15cm}}
\caption{#2}
  \label{fig:#3}
\end{figure}
\fi
}

\newcommand{\diag}{{\rm diag}}

\usepackage{color}
\usepackage{lipsum, enumitem}% http://ctan.org/pkg/lipsum
\usepackage{float}% http://ctan.org/pkg/float
\newcommand\blfootnote[1]{%
  \begingroup
  \renewcommand\thefootnote{}\footnote{#1}%
  \addtocounter{footnote}{-1}%
  \endgroup
}

\def\bu{ \textbf{u} }

\def\ba{ \textbf{a} }

\def\bC{ {\mathbf{C}} }

\def\bH{ {\mathbf{H}} }
\def\bI{ {\mathbf{I}} }

\def\bX{ {\mathbf{X}} }

\def\blambda{ \bm{\lambda} }

\def\bLambda{ {\mathbf{\Lambda}} }

\def\mcalN{ \mathcal{N} }
\def\mcalT{ \mathcal{T} }
\def\nn{{ \parallel   }}
\def\RR{{ \mathbb{R}  }}

\def\NN{{ \mathbb{N}  }}

\def\diag{{ \text{diag}   }}

\def\supp{{ \mathbf{supp}  }}
\def\bd{{ \mathbf{d}  }}
\def\bf{{ \mathbf{f}  }}
\def\be{{ \mathbf{e}  }}
\def\bx{{ \mathbf{x}  }}

\def\bz{{ \mathbf{z}  }}

\def\bg{ \mathbf{g}  }

\newtheorem{proposition}{Proposition}

\begin{document}
\title{Distributed Optimization of Multi-Beam Directional Communication Networks}

\author{\IEEEauthorblockN{Theodoros Tsiligkaridis}
\IEEEauthorblockA{MIT Lincoln Laboratory\\
Lexington, MA 02141, USA\\
Email: \texttt{ttsili@ll.mit.edu}}}


\maketitle

\blfootnote{This work is sponsored by Air Force Contract FA8721-05-C-0002. Opinions, interpretations, conclusions and recommendations are those of the author and are not necessarily endorsed by the Department of Defense or the U.S. Government.

Approved for Public Release, Distribution Unlimited.
}

\begin{abstract}
We formulate an optimization problem for maximizing the data rate of a common message transmitted from nodes within an airborne network broadcast to a central station receiver while maintaining a set of intra-network rate demands. Assuming that the network has full-duplex links with multi-beam directional capability, we obtain a convex multi-commodity flow problem and use a distributed augmented Lagrangian algorithm to solve for the optimal flows associated with each beam in the network. For each augmented Lagrangian iteration, we propose a scaled gradient projection method to minimize the local Lagrangian function that incorporates the local topology of each node in the network. Simulation results show fast convergence of the algorithm in comparison to simple distributed primal dual methods and highlight performance gains over standard minimum distance-based routing.
\end{abstract}


% \fi
\IEEEpeerreviewmaketitle

\section{Introduction}
Missions where multiple communication goals are of interest are becoming more prevalent in military applications. Multilayer communications may occur within a coalition; for example, a team consisting of ground vehicles and an airborne set of assets may desire to maximize data rate for communication from the ground team to the airborne team and vice-versa while simultaneously maintaining various types of communications within each team. This paper considers such a particular scenario where a team of nodes capable of multi-beam directional communications wants to send a common message to another team, while maintaining various routing demands within pairs of nodes. The optimization variables are the powers allocated to each beam and correspond to flows across edges in the network. %Previous work on combining communication and jamming includes time-division/frequency-division approaches which provides limited throughput and jamming gains, transmit beamforming with superimposed artificial noise which requires enough degrees of freedom to be effective and may not use channel knowledge optimally. Recent work on joint communication and jamming over the same aperture and same band that optimize the tradeoff between the two objectives has been shown to be effective in single-user and multi-user settings \cite{SCEW:2016}. The multiuser work has focused on pair-to-pair links so far.

Digital antenna arrays capable of adaptive multi-beam communications are getting closer to practice. Using transmit and receive beamforming, these arrays can form multiple links over the same frequency and may even steer nulls to mitigate interference from nearby nodes yielding higher spatial reuse and enhancing network throughput \cite{Litva:1996}. These advanced PHY techniques have to be coupled with a MAC layer. Recent prior work on this networking topic considered different MAC layer policies and proposed a new uncoordinated random access MAC policy for such systems and evaluated its throughput performance as a function of various parameters \cite{Kuperman:2016}. The recent work of \cite{MacLeod:2016} explored various spatial processing strategies for multi-beam transmission with a simple MAC layer in a simulation study. Full-duplex communications ease the burden on the MAC layer design and are projected to soon become a reality \cite{Johnston:2013, Han:2014} and can significantly increase network capacity \cite{Wang:2016} at the expense of carefully managing self-interference \cite{Quan:2017}.

In this paper, we assume full-duplex capabilities and ideal multi-beam technology in order to focus on higher level issues of distributed optimization of resources for several tasks. We adopt a constrained optimization framework to balance the power tradeoff for maximizing data rate for a common message sent to a central station and maintaining data demands for routing packets within the network. The optimal flows that arise from the solution of the optimization problem can be used to route different messages across the network. We simplify this optimization problem into a convex multicommodity flow problem and solve using a distributed augmented Lagrangian (AL) decomposition technique. We derive a scaled gradient algorithm to minimize each local AL function at each node and rigorously show its implementation requires local computations. Simulations are performed on an illustrative example to highlight gains over standard minimum distance routing protocols, and also show convergence rate improvements over simple primal dual optimization methods.

\section{Problem Formulation} \label{sec:model}
The blue team network is composed on $N$ nodes, and is modeled as a graph $\mathcal{G} = (\mcalN,\mathcal{E})$. The graph may be undirected or directed. Node $i$ may transmit data to node $j$ if $(i,j)\in \mathcal{E}$. We let $R_m$ denote a desired data rate for the traffic originating at a source node $i_m$ and ending at a sink node $j_m$. Each of these data traffic demands must be satisfied within the blue team. Define $P_{ij}$ as the power allocated by node $i$ for transmitting data to node $j$, and its associated channel capacity $c_{ij}(P_{ij}) = \log_2(1+f_{ij} P_{ij})$ bits/s/Hz for a path loss constant $f_{ij}$. Let $P_{\text{max}}$ denote the maximum transmit power of a particular node. Each node $i$ may allocated power for communications within its own team equal to $\sum_{j:(i,j)\in \mathcal{E}} P_{ij}$, and power for communication to the central station, $P_{i,C} = P_{max} - \sum_{j:(i,j)\in \mathcal{E}} P_{ij}$. We define $x_{ij}(m)$ as the information flow from node $i$ to $j$ for commodity $m$, where the total number of commodities are $M$. Each commodity here corresponds to a desired message to be sent from one node to the other.

The links are line-of-sight with path loss $f_{ij} = \frac{1}{N_0 W (4\pi/\lambda)^2 d_{ij}^2}$, where $d_{ij}$ is the distance between nodes $i$ and $j$, and $N_0 = kT$ is the noise figure of the receiver. The path loss from node $i$ to the central station is denoted as $f_{i,C}$ and is similarly defined.

We seek to optimally allocate power among several transmit beams per node in order to maximize the total signal-to-interference noise ratio at the central station receiver since the data rate is given by $R_{C} = W \log_2(1 + \sum_{i\in \mcalN} P_{i,C} f_{i,C})$ bits/s. This framework may be generalized further to include multiple central stations. The optimization problem is as follows:
\begin{align*}
	& \max_{\{x_{ij}(m), P_{ij}, P_{i,C} \}} \sum_{i\in \mcalN} P_{i,C} f_{i,C} \\
	&\text{subject to }  \sum_{m=1}^M x_{ij}(m) \leq c_{ij}(P_{ij}), \quad \forall (i,j) \in \mathcal{E} \\
		& \sum_{j: (i,j) \in \mathcal{E}} x_{ij}(m) - \sum_{j: (j,i) \in \mathcal{E}} x_{ji}(m) =
			\begin{cases}
				R_m, & \text{if } i = i_m \\
				-R_m, & \text{if } i = j_m \\
				0,              & \text{otherwise}
			\end{cases}, \\
		&\qquad \qquad  \qquad  \qquad  \qquad \forall i\in \mcalN, \forall m\in \mathcal{M} \\
		& \sum_{j: (i,j) \in \mathcal{E}} P_{ij} + P_i^J \leq P_{max}, \quad \forall i\in \mcalN \\
		& x_{ij}(m) \geq 0, P_{ij} \geq 0, P_i^J \geq 0
\end{align*}
This problem is equivalent to the primal optimization problem shown below by eliminating variables $\{P_{i,C}\}$.
\begin{align*}
	& \textbf{(P)} \qquad  \min_{\{x_{ij}(m), P_{ij} \}} \sum_{(i,j)\in \mathcal{E}} P_{ij} f_{i,C}  \\
	& \text{subject to } \sum_{m=1}^M x_{ij}(m) - c_{ij}(P_{ij}) \leq 0 \\
		& \sum_{j: (j,i) \in \mathcal{E}} x_{ji}(m) - \sum_{j: (i,j) \in \mathcal{E}} x_{ij}(m) + s_i(m) = 0  \\
		& \sum_{j: (i,j) \in \mathcal{E}} P_{ij} - P_{max} \leq 0  \\
		& x_{ij}(m) \geq 0, P_{ij} \geq 0
\end{align*}
where
\begin{equation*}
	s_i(m) \defequal \begin{cases}
									R_m, & \text{if } i = i_m \\
									-R_m, & \text{if } i = j_m \\
									0,              & \text{otherwise}
								\end{cases}
\end{equation*}
We note that the primal problem \textbf{(P)} is a convex optimization problem. For large enough $P_{max}$ and appropriate rates $R_m$, Slater's condition can be shown to hold. This implies strong duality and existence of optimal dual and primal solutions. We note however that there is no simple way to test if a solution exists to this problem for an arbitrary $P_{max}$ and desired rates $R_m$.

Next, we reformulate the problem \textbf{(P)} by eliminating the variables $\{P_{ij}\}$ by making use of the following proposition.
\begin{proposition} \label{prop:eliminateP}
	Assume that an optimal solution to the convex optimization problem \textbf{(P)} exists and is given by $(X^*,P^*)$. Then, we must have 
	\begin{equation} \label{eq:nec}
		\sum_{m=1}^M x_{ij}^*(m) = c_{ij}(P_{ij}^*), \forall (i,j) \in \mathcal{E}
	\end{equation}
\end{proposition}
\begin{proof}
	Assume that (\ref{eq:nec}) is violated. Then there exists an edge $(i',j')\in \mathcal{E}$ and $\epsilon > 0$ such that
	\begin{equation*}
		\sum_{m=1}^M x_{i'j'}^*(m) = c_{i'j'}(P_{i'j'}^*) - \epsilon
	\end{equation*}
	Now, set 
	\begin{equation*}
		\tilde{P}_{i',j'} \defequal c_{i',j'}^{-1}\left( \sum_m x_{i'j'}^*(m) \right)
	\end{equation*}
	By the continuity and monotonicity of the function $c_{i'j'}(\cdot)$, $\tilde{P}_{i'j'} < P_{i'j'}^*$. Define the new solution $(X^*,\tilde{P}=\{\tilde{P}_{i'j'},\{P_{ij}^*\}_{(i,j)\neq (i',j')}\})$. This solution also satisfies all the constraints, but achieves a better objective function value, i.e., $\tilde{P}_{i'j'} f_{i',C} + \sum_{(i,j)\neq (i',j')} P_{ij}^* f_{i,C} < \sum_{(i,j)\in \mathcal{E}} P_{ij}^* f_{i,C}$. This contradicts the optimality of $(X^*,P^*)$. The proof is complete.
\end{proof}

Proposition \ref{prop:eliminateP} allows us to eliminate the power variables $P_{ij}$ and transform the problem \textbf{(P)} into a multi-commodity optimization problem:
\begin{empheq}[box=\fbox]{align*}
   & \textbf{(Q)} \qquad \min_{\{x_{ij}(m)\}} \sum_{(i,j)\in \mathcal{E}} w_{ij} \exp\left( \ln(2) \sum_{m=1}^M x_{ij}(m) \right) \\
	 & \text{subject to } \sum_{j: (j,i) \in \mathcal{E}} x_{ji}(m) - \sum_{j: (i,j) \in \mathcal{E}} x_{ij}(m) + s_i(m) = 0 \\
		& x_{ij}(m) \geq 0
\end{empheq}
where $w_{ij} \defequal \frac{f_{i,C}}{f_{ij}}$ are positive weights.

The problem \textbf{(Q)} is also a convex optimization problem since the constraints are convex and the objective is a sum of convex functions, each one being a composition of a convex function ($\phi(u)=e^u$) with a linear combination of the variables ($y_{ij} = \sum_{m} x_{ij}(m)$). This optimization can be interpreted as trying to minimize the flows along each arc based on relative importance weights $f_{i,C}/f_{ij}$ subject to the flow conservation constraints. For large $f_{i,C}/f_{ij}$, the flow along arc $(i,j)$ tends to be minimized in order to put more power into communicating to the central station since node $i$ tends is closer to the station than the node $j$ it is communicating with, while for small $f_{i,C}/f_{ij}$ there is less emphasis placed on communicating to the station since node $i$ is closer to node $j$. The relative path loss ratio $f_{i,C}/f_{ij}$ controls the tradeoff between communication to the central station and node $j$ for each node pair $(i,j)$.

Once the optimal solution to \textbf{(Q)} is obtained, the optimal power levels may be obtained using:
\begin{equation*}
	P_{ij}^* \defequal c_{ij}^{-1}\left( \sum_m x_{ij}^*(m) \right) = \frac{e^{\ln(2) \sum_m x_{ij}^*(m)} - 1}{f_{ij}}
\end{equation*}
and feasibility is easily checked by ensuring $\sum_{j: (i,j)\in \mathcal{E}} P_{ij}^* \leq P_{max}$ for all nodes $i\in \mcalN$. This solution coincides with the solution of \textbf{(P)}.

\section{Simple Distributed Primal-Dual Algorithm for Solving \textbf{(Q)}} \label{sec:primaldualX}
The Lagrangian for problem \textbf{(Q)} is given by:
\begin{align*}
	&L(X,p) = \sum_{(i,j) \in \mathcal{E}} w_{ij} \exp\left( \ln(2) \sum_m x_{ij}(m) \right) \\
		&+ \sum_m \sum_{i\in \mcalN} p_i(m) \left( \sum_{j: (j,i)\in \mathcal{E}} x_{ji}(m) - \sum_{j: (i,j)\in \mathcal{E}} x_{ij}(m) + s_i(m) \right)
\end{align*}
where $p$ are the Lagrange multipliers associated with the flow of conservation. Using the primal-dual method approach of \cite{Nedic:1:2009, Nedic:2:2009} to find approximate solutions to \textbf{(Q)}, the primal and dual updates become:
%\vspace{3mm}
\noindent \textbf{Primal Variable Updates}
\begin{align*}
	x_{ij}(m)^{k+1} &= \left[ x_{ij}(m)^k - \alpha \frac{\partial L(X,p)}{\partial x_{ij}(m)} \right]_+ \\
		& = \Bigg[ x_{ij}(m)^k - \alpha \Big( \ln(2) w_{ij} e^{\ln(2) \sum_{m'} x_{ij}(m')^k} \\
		&\qquad \qquad + p_j(m)^k - p_i(m)^k \Big) \Bigg]_+
\end{align*}
\noindent \textbf{Dual Variable Updates}
\begin{align*}
	&p_i(m)^{k+1} = p_i(m)^k + \alpha \frac{\partial L(X,p)}{\partial p_i(m)} \\
		&= p_i(m)^k + \alpha \left( \sum_{j: (j,i) \in \mathcal{E}} x_{ji}(m)^k - \sum_{j: (i,j) \in \mathcal{E}} x_{ij}(m)^k + s_i(m) \right)
\end{align*}
where $\alpha>0$ is a small step size.

The primal solution after $k$ iterations is obtained by averaging the iterates:
\begin{equation*}
	\hat{x}_{ij}(m)^{(k)} = \frac{1}{k} x_{ij}(m)^k + \left( 1 - \frac{1}{k} \right) \hat{x}_{ij}(m)^{(k-1)}
\end{equation*}
It is known that this primal solution will be feasible asymptotically. We note that the averages may be implemented recursively as shown above to save memory. Under some mild conditions, this iterative algorithm converges to a saddle point of the Lagrangian function. We remark that these updates can be computed using local computations, so the algorithm may be implemented in a decentralized manner.

\section{Distributed Augmented Lagrangian Method for Solving \textbf{(Q)}} \label{sec:aug_lagrangian}
The objective function of problem \textbf{(Q)} is convex, twice-differentiable, and monotonically nondecreasing with respect to the flows $\{ x_{ij}(m) \}$. Furthermore, the only constraints present are flow conservation constraints and the flow nonnegativity constraints. This makes the path flow formulation applicable \cite{Bertsekas:Network:1998}. From the monotonicity of the objective function, an optimal flow vector may be constructed using only simple path flows. Optimizing over the path flows is not very practical since the paths from a source to a destination need to be enumerated first. The problem of finding all such simple paths cannot be accomplished in polynomial time. Furthermore, these types of algorithm in addition to the path augmentation, blocking flows and linear programming require global coordination.

%Instead, we use duality to find the optimal multipliers using unconstrained optimization and then use that to find the optimal flows in a distributed manner.
We adopt the augmented Lagrangian (AL) algorithm of \cite{Chatzipanagiotis:2015}. Define $\bx_i(m)=[x_{i1}(m),\dots,x_{iN}(m)]^T$ as the vector of flows of commodity $m$ that node $i$ routes towards all other nodes $j$, and $\bx_i = [\bx_i(1)^T,\dots,\bx_i(M)^T]^T$ is the collection of all such vectors making up a local variable of node $i$. The demand vector for commodity $m$ is given by $\bd(m)=[d_1(m),\dots,d_N(m)]^T$. This vector is defined as $d_i(m)=+R_m$ for $i=i_m$, $d_i(m)=-R_m$ for $i=j_m$ and $d_i(m)=0$ for the remaining nodes $i$. Define the local neighborhood of node $i$ as $\mcalN_i$ (not including $i$), and the two-hop neighborhood of node $i$ as $\mcalT_i$. Note that $\mcalN_i \subseteq \mcalT_i$.

Problem \textbf{(Q)} can be equivalently written as:
\begin{align*}
	& \textbf{(Q')} \qquad \min_{\{x_{ij}(m)\}} \sum_{i\in \mcalN} \sum_{j: (i,j)\in \mathcal{E}} w_{ij} 2^{\sum_{m=1}^M x_{ij}(m)}  \\
	&\text{subject to } \sum_{i \in \mcalN} \bC_i \bx_i(m) = \bd(m), \forall m\in \mathcal{M} \\
		&\qquad \bx_i(m) \succeq 0, \forall i \in \mathcal{N}, \forall m\in \mathcal{M}
\end{align*}
where the matrix $\bC_i \in \RR^{N\times N}$ is defined as
\begin{equation*}
	[\bC_i]_{l,k} = \begin{cases} +1, & l=i, k\in \mathcal{N}_i \\ -1, & l=k,k\in \mathcal{N}_i \\ 0, & \text{else} \end{cases}
\end{equation*}
In matrix form, $\bC_i = \bf_i \otimes \be_i - \diag(\bf_i)$, where $[\bf_i]_{\mathcal{N}_i} = 1$ and zero elsewhere.

The local AL of node $i$ at iteration $k$ is given by:
\begin{align}
	&\Lambda_i(\bx_i; \{\bx_j^{k}\}_{j \in \mathcal{T}_i}, \{\lambda_j^{k}(m)\}_{j\in \mathcal{N}_i \cup \{i\}}) \nonumber \\
		&\quad = \sum_{j: (i,j)\in \mathcal{E}} w_{ij} 2^{\ba_j^T \bx_i} + \sum_m (\lambda(m)^k)^T \bC_i \bx_i(m) \nonumber \\
		&\quad + \frac{\rho}{2} \sum_m \nn \bC_i\bx_i(m) + \sum_{j\neq i} \bC_j \bx_j(m)^k - \bd(m) \nn_2^2 \label{eq:localAL}
\end{align}
where $\rho>0$ is a regularization parameter. The selection vectors are given as $\ba_j = \mathbf{1}_M \otimes \be_j$. The accelerated distributed augmented Lagrangian (ADAL) method of \cite{Chatzipanagiotis:2015} is summarized in Algorithm \ref{alg:adal}.
\begin{algorithm}[tb]
   \caption{Accelerated Distributed Augmented Lagrangian (ADAL) Algorithm}    \label{alg:adal}
\begin{algorithmic}
		\STATE {\bfseries Input:} initialize $\lambda$ and flow variables $\bX$.
		\STATE Set $\rho > 0$ and $\tau \in (0,1/d_{max})$, where $d_{max}$ is the maximum degree in the graph.
		\STATE 1. Each node $i$ computes the minimizer:
			\begin{equation} \label{eq:minLocalAL}
				\hat{\bx}_i^k = \arg \min_{\bx_i \succeq 0} \Lambda_i(\bx_i; \{\bx_j^k\}_{j \in \mathcal{T}_i}, \{\lambda_j^k(m)\}_{j\in \mathcal{N}_i})
			\end{equation}
		\STATE 2. For every node $i$, update local flow vector:
			\begin{equation*}
				\bx_i^{k+1} = \bx_i^k + \tau (\hat{\bx}_i^k - \bx_i^k)
			\end{equation*}
		\STATE 3. Update Lagrange multiplier for each node $i$ and commodity $m$:
			\begin{align*}
				&\lambda_i^{k+1}(m) = \lambda_i^k(m) \\
					&\quad + \rho \tau \left( \sum_{j: (i,j) \in \mathcal{E}} x_{ij}^{k+1}(m) - \sum_{j: (j,i) \in \mathcal{E}} x_{ji}^{k+1}(m) - d_i(m) \right)
			\end{align*}
			Go to Step 1.
\end{algorithmic}
\end{algorithm}

\subsection{Minimizing Local AL}
The minimization problem (\ref{eq:minLocalAL}) can be solved by a projected gradient method \cite{Bertsekas:NonlinearProg:1998}. We assume an initial condition of $\bx_i(0) = \bx_i^k$ and perform iterative updates on the flow vector using Algorithm \ref{alg:scaled_proj_grad}.
\begin{algorithm}[tb]
   \caption{Scaled Projected Gradient Algorithm for solving (\ref{eq:minLocalAL}) }    \label{alg:scaled_proj_grad}
\begin{algorithmic}
		\STATE Set $s > 0$ and $\beta,\sigma \in (0,1)$.
		\REPEAT
		\STATE For $t\geq 1$ do:
			\begin{align*}
				\bar{\bx}_i(t) &= \left[ \bx_i(t) + s \bd(t) \right]_+ \\
				\bx(t+1) &= \bx(t) + \alpha(t) (\bar{\bx}(t) - \bx(t))
			\end{align*}
			\STATE where $\bu(t) = \bar{\bx}_i(t) - \bx_i(t)$, $\bd(t)$ is the scaled gradient direction evaluated at $\bx_i(t)$, $\alpha(t) = \beta^{m(t)}$ and
			\begin{align*}
				m(t) = \min\Big\{ m\in \NN: \Lambda_i(\bx_i(t)) - &\Lambda_i(\bx_i(t) + \beta^m \bu(t)) \\
					&\geq \sigma \beta^m \left< \bd(t),\bu(t) \right> \Big\}
			\end{align*}
		\UNTIL{Convergence}
\end{algorithmic}
\end{algorithm}
In practice to keep computational complexity bounded, we perform only several projected gradient descent iterations until a local convergence criterion is satisfied, i.e., $\nn [\bx_i^{(k)} - \nabla \Lambda_i(\bx_i^{(k)})]_+ - \bx_i^{(k)} \nn_2 \leq \epsilon$ for some small $\epsilon$.

\subsubsection{Gradient Calculation}
The gradient takes the form:
\begin{align*}
	\nabla_{\bx_i} \Lambda_i(\bx_i) &= \sum_{j \in \mcalN_i} w_{ij} 2^{\ba_j^T \bx_i} \ln(2) \ba_j \\
		&\quad + (\bI_M \otimes \bC_i^T) (\blambda^k + \rho \bz_i^k) + \rho (\bI_M \otimes \bC_i^T\bC_i) \bx_i
\end{align*}
where $\blambda^k = [\blambda(1)^T,\dots,\blambda(M)^T]^T \in \RR^{NM}$, $\bz_i^k = [(\bz_i(1)^k)^T,\dots,(\bz_i(M)^k)^T]^T \in \RR^{NM}$ and $\bz_i(m)^k = \sum_{j \neq i} \bC_j \bx_j^k - \bd(m) \in \RR^N$.

% proof for needed only 2-hop neighbors to compute gradient
We next prove that the gradient may be calculated using local information within two-hop neighbors $j \in \mathcal{T}_i$. The gradient vector may be decomposed as:
\begin{equation*}
	\nabla_{\bx_i} \Lambda_i(\bx_i) = \begin{bmatrix} \nabla_{\bx_i(1)} \Lambda_i(\bx_i) \\ \vdots \\ \nabla_{\bx_i(M)} \Lambda_i(\bx_i) \end{bmatrix}
\end{equation*}
where
\begin{align*}
	\nabla_{\bx_i(m)} \Lambda_i(\bx_i) &= \sum_{j \in \mcalN_i} w_{ij} 2^{\ba_j^T \bx_i} \ln(2) \be_j + \bC_i^T (\blambda(m) + \rho \bz_i(m)^k) \\
		&\qquad + \rho \bC_i^T\bC_i \bx_i
\end{align*}
Since the rows $[\bC_i^T]_{l,\cdot} = 0$ for all $l\notin \mcalN_i$, it follows that $\supp(\nabla_{\bx_i(m)} \Lambda_i(\bx_i)) = \mathcal{N}_i$.

We first note that the term $\bC_i^T\blambda(m)$ is locally computable since:
\begin{equation} \label{eq:l1}
	[\bC_i^T\blambda(m)]_l = \begin{cases} \lambda_i(m) - \lambda_l(m), & l\in \mcalN_i \\ 0, & \text{else}  \end{cases}
\end{equation}
so access to $\lambda_l(m)$ from all neighbors $l\in \mcalN_i$ is only needed.

We then focus our attention to the term $\bC_i^T\bC_i \bx_i$. This is locally computable since:
\begin{equation*}
	[\bC_i^T\bC_i \bx_i]_l = \begin{cases} 2x_{il} + \sum_{j\in \mcalN_i\backslash \{l\}} x_{ij}, & l\in \mcalN_i \\ 0, & \text{else}  \end{cases}
\end{equation*}

Finally, we consider the term $\bC_i^T \bz_i(m)^k$, where $\bz_i(m)^k = \sum_{j\neq i} \bC_j \bx_j^k - \bd(m)$. Since $\supp(\bC_i^T \bz_i(m)^k) = \mcalN_i$ from (\ref{eq:l1}), we only focus on calculating this vector for coordinates $l \in \mcalN_i$. Then, for $l\in \mcalN_i$, we have:
\begin{align*}
	&[\bC_i^T \bz_i(m)^k]_l = [\bz_i(m)^k]_i - [\bz_i(m)^k]_l \\
		&\quad = \left( \sum_{j \neq i} [\bC_j\bx_j^k]_i - \sum_{j \neq i} [\bC_j\bx_j^k]_l \right) - \left( d_i(m) - d_l(m) \right)
\end{align*}
Next, we show that $\sum_{j \neq i} [\bC_j\bx_j]_i - \sum_{j \neq i} [\bC_j\bx_j]_l$ involves only $x_{j,t}$ for $j\in \mcalT_i$.
\begin{align*}
	&\sum_{j \neq i} [\bC_j\bx_j]_i - \sum_{j \neq i} [\bC_j\bx_j]_l = \sum_{j \in \mcalN_i} [\bC_j\bx_j]_i - \sum_{j \in \mcalT_i} [\bC_j\bx_j]_l \\
		&\quad = \sum_{j\in \mcalN_i} ([\bC_j\bx_j]_i-[\bC_j\bx_j]_l) - \sum_{j \in \mcalT_i\backslash \mcalN_i} [\bC_j\bx_j]_l \\
		&\quad = \sum_{j\in \mcalN_i} \left(-x_{ji} - \begin{cases} \sum_{t\in \mcalN_j} x_{jt}, & l=j \\ -x_{jl}, & l\neq j \end{cases} \right) + \sum_{j\in \mcalT_i\backslash \mcalN_i} x_{jl}
\end{align*}
Thus, to compute the gradient, nodes only need information from their two-hop neighbors and do not require global knowledge of the network.

\subsubsection{Hessian Calculation}
The Hessian matrix of $\bLambda_i$ takes the form:
\begin{equation*}
	\nabla_{\bx_i}^2 \bLambda_i(\bx_i) = \sum_{j\in \mathcal{N}_i} w_{ij} 2^{\ba_j^T\bx_i} (\ln(2))^2 \ba_j \ba_j^T + \rho (\bI_M \otimes \bC_i^T \bC_i)
\end{equation*}
The Hessian is rank-deficient with rank $M|\mathcal{N}_i|$, depending on the size of the local neighborhood. The full $MN\times MN$ Hessian matrix can be decomposed as:
\begin{equation*}
	\nabla_{\bx_i}^2 \bLambda_i(\bx_i) = \begin{bmatrix} \nabla_{\bx_i(1)}^2 \bLambda_i(\bx_i) & \cdots & \nabla_{\bx_i(1)} \nabla_{\bx_i(M)}^T \bLambda_i(\bx_i) \\ \vdots & \ddots & \vdots \\ \nabla_{\bx_i(M)} \nabla_{\bx_i(1)}^T \bLambda_i(\bx_i) & \cdots & \nabla_{\bx_i(M)}^2 \bLambda_i(\bx_i) \end{bmatrix}
\end{equation*}
where each blockwise component is
\begin{equation*}
	\nabla_{\bx_i(m)} \nabla_{\bx_i(m')}^T \bLambda_i(\bx_i) = \sum_{j\in \mathcal{N}_i} w_{ij} 2^{\ba_j^T\bx_i} (\ln(2))^2 \be_j \be_j^T + \rho \bC_i^T\bC_i I_{\{m=m'\}}
\end{equation*}
The matrix $\bC_i^T\bC_i$ depends on the local neighborhood $\mathcal{N}_i$ as:
\begin{equation*}
	[\bC_i^T \bC_i]_{l,k} = \begin{cases} 2, & l=k, k\in \mathcal{N}_i \\ 1, & l\neq k,l\in \mathcal{N}_i,k\in \mathcal{N}_i \\ 0, & \text{else} \end{cases}
\end{equation*}
We also note that the reduced-dimension Hessian (restricted to $\mathcal{N}_i$) is always positive definite since $[\bC_i^T\bC_i]_{\mathcal{N}_i}  = \mathbf{1}_{|\mathcal{N}_i|} \mathbf{1}_{|\mathcal{N}_i|}^T + \bI_{|\mathcal{N}_i|}$. Each submatrix has support: $\supp(\nabla_{\bx_i(m)} \nabla_{\bx_i(m')}^T \bLambda_i(\bx_i)) = \mathcal{N}_i \times \mathcal{N}_i$.

\subsubsection{Scaled Gradient Direction Calculation}
Let $\mathcal{I}_i$ denote the set of indexes corresponding to the neighborhood of node $i$, i.e., $\mathcal{I}_i = \{\mathcal{N}_i+(m-1)N: m=0,\dots,M-1\}$. The reduced-dimension Hessian $\bH_{\mathcal{I}_i,\mathcal{I}_i} = [\nabla^2 \bLambda_i]_{k,l\in \mathcal{I}_i}$ has full-rank, and the reduced-dimension gradient $\bg_{\mathcal{I}_i} = [\nabla \bLambda_i]_{k \in \mathcal{I}_i}$ has nonzero components in general. Using this index set, we may work over a reduced-dimension space and calculate a scaled gradient-descent direction as:
\begin{equation*}
	\tilde{\bd} = - (\diag(\bH_{\mathcal{I}_i,\mathcal{I}_i}))^{-1} \bg_{\mathcal{I}_i} = -\frac{1}{2\rho} \bg_{\mathcal{I}_i}
\end{equation*}
with its expanded version defined as $d_l = \tilde{d}_k I_{\{l \in \mathcal{I}_i\}}$ where $l=k+(m-1)N$ for appropriate $k\in \mathcal{N}_i$ and $m\in\{0,\dots,M-1\}$. The diagonal approximation to the Hessian becomes a simple scaling of the gradient. Using this tends to have faster convergence than the unscaled gradient projection method, as the experiments show.

\section{Simulation Results}
We simulate a two-dimensional scenario where a network broadcasts a common message to a central station while maintaining desired data rates among two origin-destination node pairs using multiple beams. A center frequency of $f_c=1$ GHz and bandwidth of $W=5$ MHz is used, with a peak transmit power per node of $P_{max}=100$ Watts. In the simulation, we compare our distributed optimization algorithm (denoted as ADAL) with the OSPF routing protocol which routes messages at a certain rate through the minimum distance route, which is efficiently obtained using Dijkstra's algorithm. 

%Furthermore, we observe that the optimized approach outperforms OSPF by more than a factor of $3$ in terms of total power consumption for communications at all distances, and outperforms OSPF by a factor ranging from $1.5$ to $3$ in terms of the interference level at the red receiver depending on its location.

\begin{figure}[ht]
	\centering
		\includegraphics[width=0.50\textwidth]{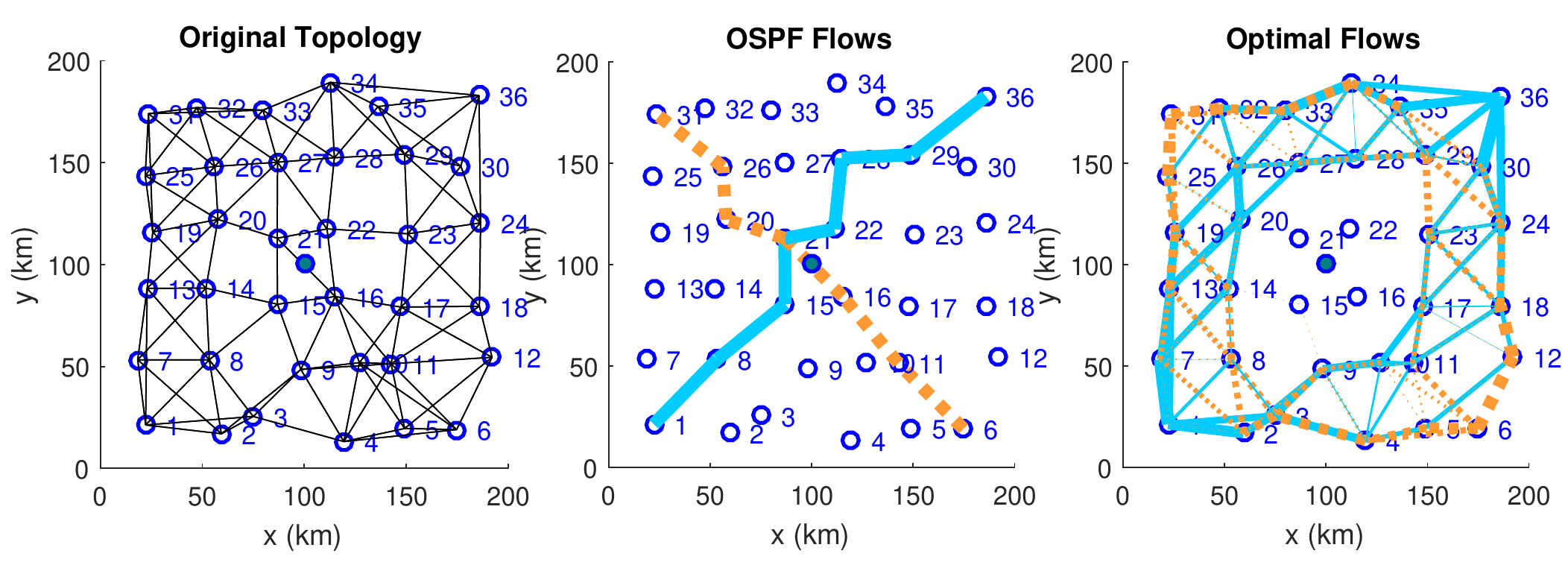}
	\caption{Approximate grid geometry. Two messages are to be sent from node $1$ to $36$ and from node $6$ to $31$, respectively, using a data rate $R=9$ bits/s/Hz. Central station receiver is located in the center of the grid. The middle panel shows the OSPF routes, and the right panel shows the optimal routes obtained using the ADAL algorithm. }
	\label{fig:grid:geom}
\end{figure}

The network consists of $N=36$ nodes arranged approximately in a $6 \times 6$ grid as the left panel of Fig. \ref{fig:grid:geom} shows. The central station is placed in the center of the grid and two messages are to be routed from node $1$ to $36$, and from node $6$ to node $31$, each with desired data rate $R=9$ bits/s/Hz. The middle panel of Fig. \ref{fig:grid:geom} shows the OSPF flows which correspond to the minimum distance routes. The right panel shows the optimal flows obtained using the ADAL algorithm, which show that routes around the central station are preferred with appropriate load balancing than full loading the shortest routes. This makes sense since the nodes in the center allocate more power for communicating to the central station and the nodes away from the center are more responsible for carrying out the intra-network communications.

\begin{figure}[ht]
	\centering
		\includegraphics[width=0.50\textwidth]{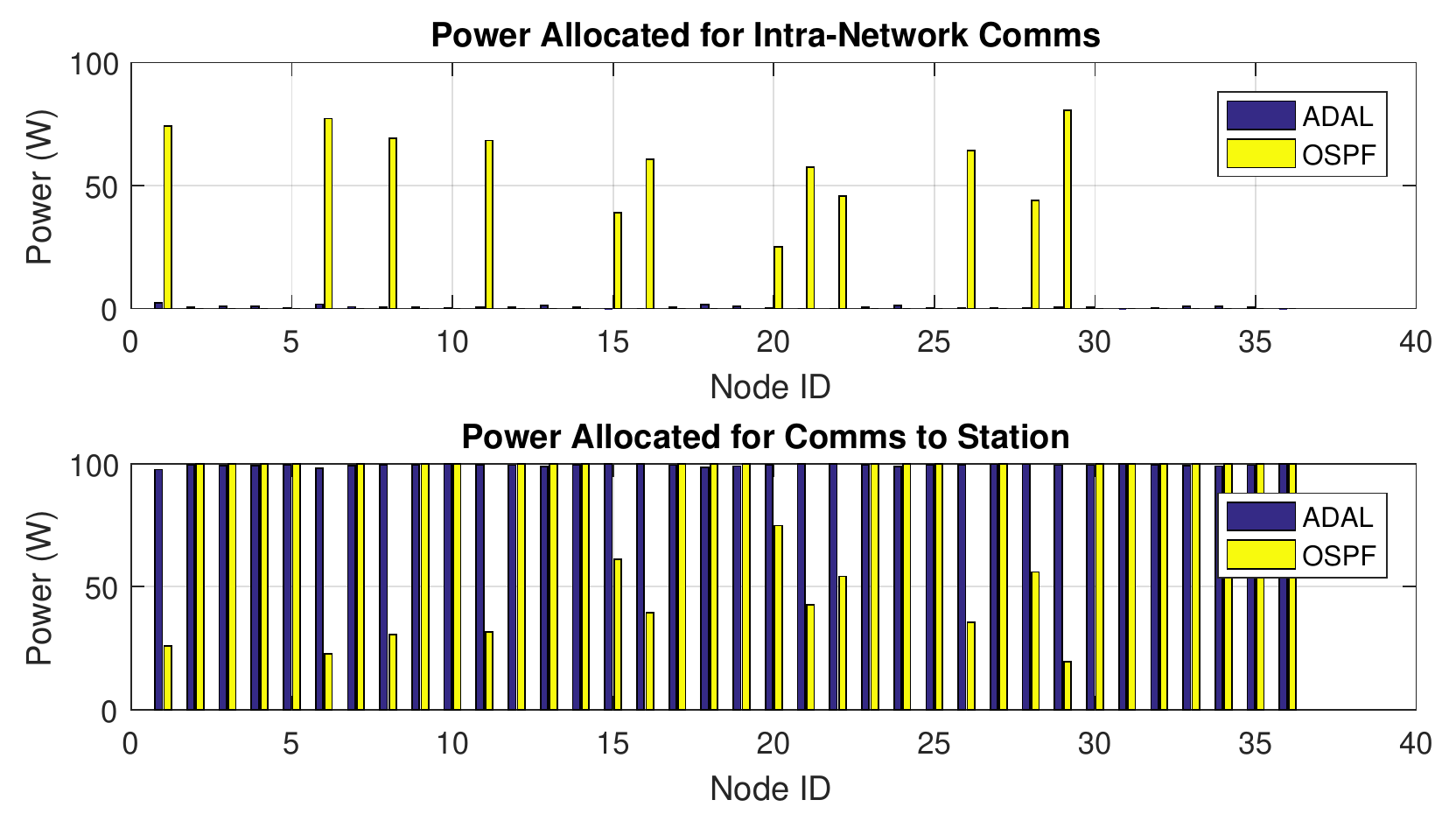}
	\caption{Power used for intra-network communications (top) and power used for communication to the central station (bottom) across nodes. ADAL uses significantly less power total for intra-network communications in comparison to OSPF. }
	\label{fig:grid:powers}
\end{figure}

The power used for intra-network communications and communications to the central station is shown in Fig. \ref{fig:grid:powers} for each node. ADAL uses significantly less power overall for communications, $22.96$W, in comparison to OSPF, $705.69$W, which is a factor of $\times 30$ reduction. As a result, nodes have more power left over to use for maximizing the data rate of the common message to the station receiver. The received power at the central station receiver for OSPF is $2.89$kW corresponding to a data rate of $R_C=68.5$Mbps, while for ADAL is $3.58$kW corresponding to a data rate of $R_C=71.9$Mbps. Depending on the distances of nodes to the central receiver and the network topology, this gap in received power may significantly boost the data rate. Our optimization algorithm performs optimal load balancing among different paths which leads to considerably less transmit power used for intra-network communications and allows for more power to be used for broadcasting the common message to the central receiver.

\begin{figure}[ht]
	\centering
		\includegraphics[width=0.50\textwidth]{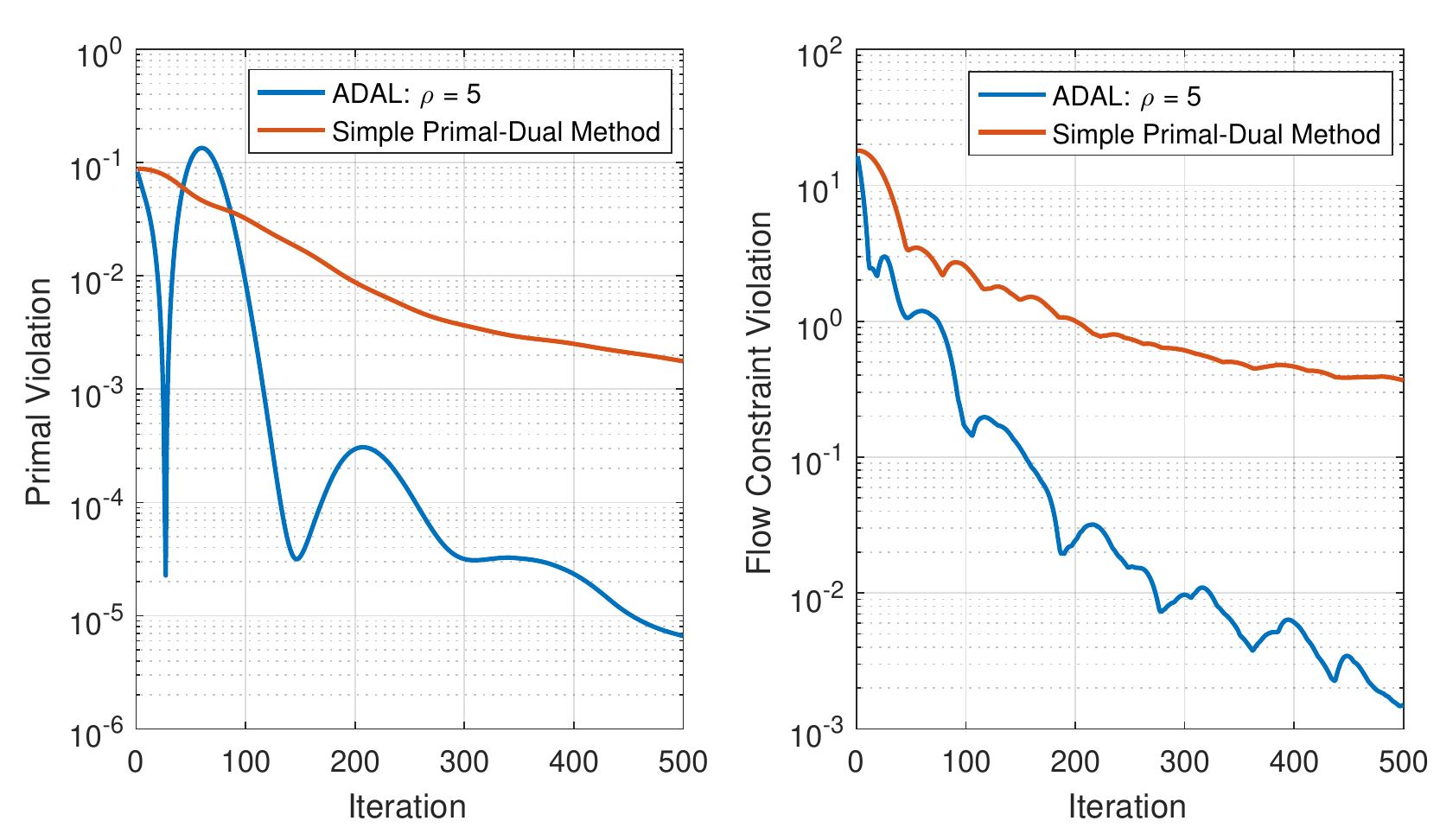}
	\caption{Convergence of objective function (left panel), and flow constraint violation (right panel). ADAL achieves faster convergence than the simple primal-dual method. }
	\label{fig:grid:convergence}
\end{figure}

Next, we examine the convergence of the augmented Lagrangian (ADAL) algorithm and compare it with the simple primal-dual method. The constraint violation metric used is
\begin{equation*}
	\sum_{m=1}^M \sum_{i \in \mcalN} \left| \sum_{j: (j,i)\in \mathcal{E}} x_{ji}^{k}(m) - \sum_{j: (i,j)\in \mathcal{E}} x_{ij}^{k}(m) + s_i(m) \right|
\end{equation*}
which is expected to converge to zero as $k$ grows. The objective function error measures the distance from the optimal primal value. Fig. \ref{fig:grid:convergence} shows that ADAL achieves significantly faster convergence than the simple primal-dual method. The simple primal-dual algorithm is fully decentralized and has guaranteed convergence to the unique minimizer but has slow convergence.

\begin{figure}[ht]
	\centering
		\includegraphics[width=0.50\textwidth]{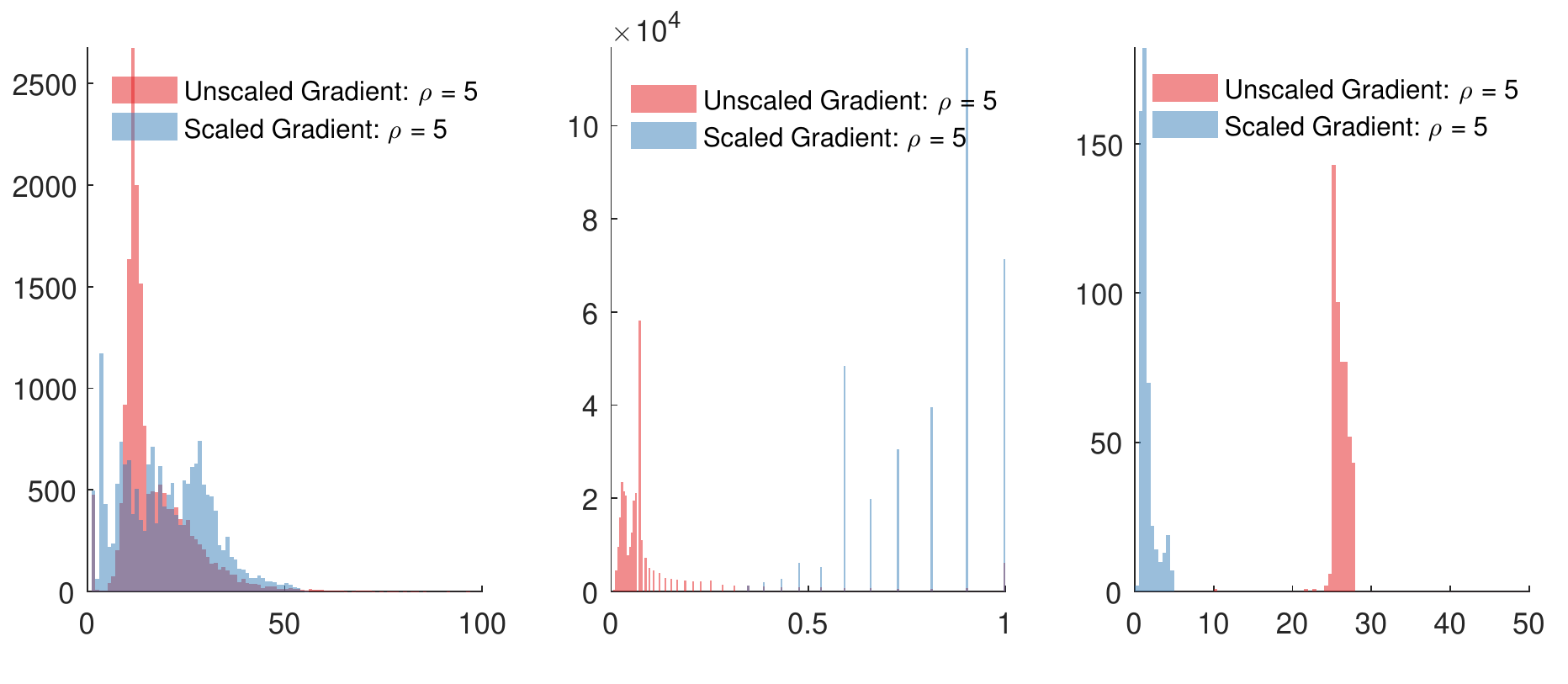}
	\caption{Histogram of number of inner iterations (left), Armijo step sizes (middle) and average Armijo steps (right) for unscaled and scaled gradient methods to minimize local AL function at each iteration of the ADAL algorithm. }
	\label{fig:grid:complexity}
\end{figure}

The computational complexity of unscaled and scaled gradient methods for minimizing the local AL function in each iteration of the ADAL algorithm is addressed next. A tolerance threshold of $\epsilon=1e^{-3}$ was set to stop the inner minimization with stopping criterion:
\begin{equation*}
	\nn [\bx_i^{(k)} - \nabla \Lambda_i(\bx_i^{(k)})]_+ - \bx_i^{(k)} \nn_2 \leq \epsilon
\end{equation*}
Figure \ref{fig:grid:complexity} on the left panel shows a histogram of the number of inner iterations needed to achieve the tolerance $\epsilon$ for the unscaled and scaled gradient methods, respectively, and they are on the same order on average. This metric determines the number of gradient evaluations of the local AL function. The middle panel shows a histogram of the number of Armijo step sizes. The scaled method tends to require step sizes close to unity while the unscaled method requires a lot of tuning and requires significantly smaller step sizes, which implies reduced latency during optimization. The right panel displays the average number of Armijo steps per inner iteration. This is an important complexity metric since it determines the number of local AL function evaluations. The scaled gradient method requires an average of only $1.5$ steps while the unscaled gradient method requires $26.1$ steps on average.

%\subsection{Ensemble Simulation}
%Next, we compare the gradient and approximate second-order versions of the ADAL algorithm with the primal-dual algorithm in terms of convergence, in addition to OSPF routing.

%\section{Simulation Results} \label{sec:simulation}
%
%%\begin{figure}[ht]
	%%\centering{
	%%\subfigure[]{
		%%\centering
		%%\includegraphics[width=0.40\textwidth]{./undirected_snapshot_iter=100.pdf}
	%%}
	%%\subfigure[]{
		%%\centering
		%%\includegraphics[width=0.40\textwidth]{./undirected_snapshot_iter=200.pdf}
	%%}
	%%\subfigure[]{
		%%\centering
		%%\includegraphics[width=0.40\textwidth]{./undirected_snapshot_iter=300.pdf}
	%%}
	%%\subfigure[]{
		%%\centering
		%%\includegraphics[width=0.40\textwidth]{./undirected_snapshot_iter=400.pdf}
	%%}
	%%}
	%%\caption{Undirected network (no jammer). Distributed control of UAV swarm for enabling relay communication from source to target nodes. Snapshots at different time instants show the evolution of the swarm. }
	%%\label{fig:snapshots_undirected}
%%\end{figure}
%%
%

\section{Conclusion}
We proposed distributed algorithms for power allocation in multibeam directional airborne networks for maximizing data rate for a common message sent by all nodes to a central station receiver while guaranteeing multiple rate demands for intra-network communications. A decomposition approach was applied to the augmented Lagrangian (AL) algorithm for solving the convex optimization problem that arises, and an efficient method for solving each subproblem in each AL iteration was presented in detail. Simulation results show the benefits of our approach in comparison to simple primal dual methods in terms of convergence. Finally, significant power savings are observed for intra-network communications with our optimized routing algorithms in comparison to standard minimum distance routing.

%\clearpage
%\vfill\pagebreak

\bibliographystyle{IEEEtran}
\bibliography{refs}

\end{document}